\documentclass[1p,nopreprintline]{elsarticle}

\usepackage{lineno,hyperref}
\modulolinenumbers[5]

\journal{ }

\usepackage{bm,graphicx}

\usepackage{amsmath,amsfonts}
\newtheorem{theorem}{Theorem}
\newtheorem{lemma}{Lemma}
\newtheorem{definition}{Definition}
\newtheorem{assumption}{Assumption}

\newtheorem{corollary}{Corollary}
\newtheorem{uc}{Use Case}

\newproof{pf}{Proof}

\begin{document}

\begin{frontmatter}

\title{Second-order cone optimization of the gradostat\tnoteref{mytitlenote}}
\tnotetext[mytitlenote]{Funding is acknowledged from the Natural Sciences and Engineering Research Council of Canada and the French LabEx NUMEV (Project ANR-10 LABX-20), incorporated into the I-Site MUSE, which partially funded the sabbatical of J.~Taylor at MISTEA lab.}

\author{Josh A. Taylor}
\ead{josh.taylor@utoronto.ca}
\address{The Edward S. Rogers Sr. Department of Electrical and Computer Engineering, University of Toronto, Toronto, Canada}

\author{Alain Rapaport}
\ead{alain.rapaport@inrae.fr}
\address{MISTEA, Univ. Montpellier, INRAE, Institut Agro, Montpellier, France}
%
%
%

\begin{abstract}
We maximize the production of biogas in a gradostat at steady state. The physical decision variables are the water, substrate, and biomass entering each tank and the flows through the interconnecting pipes. Our main technical focus is the nonconvex constraint describing microbial growth. We formulate a relaxation and prove that it is exact when the gradostat is outflow connected, its system matrix is irreducible, and the growth rate satisfies a simple condition. The relaxation has second-order cone representations for the Monod and Contois growth rates. We extend the steady state models to the case of multiple time periods by replacing the derivatives with numerical approximations instead of setting them to zero. The resulting optimizations are second-order cone programs, which can be solved at large scales using standard industrial software.
\end{abstract}


\end{frontmatter}


\section{Introduction}
The gradostat is a nonlinear dynamical system in which multiple chemostats are interconnected by mass flows and diffusion. In each chemostat, microbial growth converts a substrate to biomass. This also produces biogas, a useful energy source. Our primary motivation for this setup is the design and operation of a network of wastewater treatment plants~\cite{shen2015overview}.

We seek to maximize the production of biogas in gradostats with two standard growth rates: Monod~\cite{monod1949growth} and Contois~\cite{contois1959kinetics}. There are several physical decision variables, including the inflows of water, substrate, and biomass at each tank, and the installation of pipes between tanks. We start with a steady state model obtained by setting the derivatives in the gradostat to zero. The resulting algebraic equations specify a feasible set, which is nonconvex due to the nonlinear growth in each tank and, in some setups, discrete and bilinear mass flows between the tanks.

Nonconvex optimization can be difficult even at small scales. Second-order cone programming (SOCP) is a tractable, convex optimization class that is often used to approximate nonconvex problems~\cite{Boyd1998SOCP}. In this paper, we construct SOC and mixed-integer (MI)SOC relaxations of the gradostat. Our main original contributions, listed below, center on the convexification of the growth rate constraint.
\begin{itemize}
\item In Section~\ref{relax}, we formulate a simple relaxation of the gradostat. The relaxation is obtained by allowing the conversion of substrate to biomass to be less than or equal to the growth kinetics. In Theorem~\ref{PRexact}, we prove that when the gradostat is outflow connected, its system matrix is irreducible, and a simple condition on the growth rate is satisfied, this relaxation is exact, which is to say that the inequality is satisfied with equality.
\item In Section~\ref{PRSOC}, we identify original SOC representations of the relaxed growth constraints. Specifically, the Contois growth constraint is exactly representable as an SOC constraint. The Monod growth constraint is SOC under a constant biomass approximation. We use the convex envelopes of~\cite{tawarmalani2001semidefinite} to obtain an SOC outer approximation of the Monod growth constraint in the general case.
\item In Section~\ref{extensions}, we give two extensions. In Section~\ref{ConvexUE}, we give simple linear underestimators of the growth constraints for when the relaxations are not exact. In Section~\ref{EMP}, instead of setting the derivatives in the gradostat to zero, we replace them with linear numerical approximations. This leads to an optimization with multiple time periods, which can accommodate transient conditions.
\end{itemize}

The end result is a family of SOCPs and MISOCPs for optimizing the gradostat in steady state. Today, SOCPs with tens of thousands of variables and constraints can be solved in seconds on a typical personal computer. MISOCPs are also reasonably tractable because, like mixed-integer linear programs, there are powerful mathematical tools for speeding up their solution~\cite{drewesMISOCP2009,Atamturk2010RC,belotti2017complete}, and they are handled by industrial solvers such as Gurobi~\cite{gurobi}. This enables us to solve each MISOCP to optimality at moderate scales, typically up to one hundred binary variables.

We now review some relevant literature. We refer the reader to~\cite{smith1995theory,harmand2017chemostat} for comprehensive coverage of the chemostat. The gradostat was originally formulated in~\cite{lovitt1981gradostat} as a single series of tanks and later generalized to a network of interconnected chemostats. There is an ongoing literature stream focusing on its nonlinear analysis~\cite{jager1987competition} and control~\cite{bayen2014minimal}. To date, there have been no applications of numerical convex optimization to the gradostat in steady state. 

From a technical viewpoint, the nearest topic to ours is the optimization of chemical process networks~\cite{biegler1997systematic}. Common features include bilinear mass flows, which we similarly linearize using disjunctive programming~\cite{grossmann2002disj,lee2003global}, and quotients of variables, which we approximate with convex envelopes~\cite{tawarmalani2001semidefinite,tawarmalani2002convex,tawarmalani2013convexification} in Section~\ref{MonodCE}. The main feature of the gradostat that is not present in chemical process networks is the microbial growth in the tanks, which is the core focus of this paper.

The design of chemical reactors has been studied for decades, see, e.g., \cite{froment1990chemical}. A handful of papers within this stream have focused on optimal design, some using the chemostat and gradostat. An early text is~\cite{aris2000optimal}, which uses dynamic programming, but not the chemostat. Several later studies model growth with the Monod equation and derive analytical expressions for parameters such as concentration, residence time, and tank volume~\cite{bischoff1966optimal,luyben1982optimal,hill1989minimum,de1996bioreactors}. References~\cite{harmand2003optimal,harmand2004optimal,harmand2005optimal} build on this approach, deriving numerical and qualitative conditions for optimizing reactors modeled as `steady-state equivalent biological systems.'

More recent studies have explicitly optimized the design and operation of interconnected tanks. The volumes of series bioreactors with Monod and Contois growth are optimized in~\cite{zambrano2014optimizing,zambrano2015optimal} for a given output substrate concentration at steady state. The operation of two chemostats in series is optimized in~\cite{bayen2019steady}, and the volumes and diffusion rate of a chemostat with a side compartment are optimized in~\cite{crespo2020analysis}. Of particular relevance are \cite{ROBLESRODRIGUEZ2018880,robles2019management}, which model interconnected wastewater treatments plants as a gradostat. The latter formulates inflow management as a model predictive control problem, which it solves using particle swarm optimization. While closely related to our perspective, these papers focus on different problem statements with specific network structures, and do not employ convex relaxations or SOCP.





The rest of the paper is organized as follows. Section~\ref{setup} covers the relevant background and states the nonconvex optimization problem. Section~\ref{relax} gives a simple relaxation and analyzes when it is exact. SOC representations of the Monod and Contois growth constraints are given in Section~\ref{PRSOC}, and the bilinear mass balance constraints are linearized in Section~\ref{disjunctive}. Section~\ref{extensions} gives linear underestimators of the growth constraint and extends the models to the case of multiple time periods. In Section~\ref{app}, we summarize the models and implement those that are SOCPs and MISOCPs in numerical examples.

\section{Setup}\label{setup}
\subsection{Definitions}\label{s:def}
We consider a gradostat with $n$ interconnected tanks, and denote the set of tanks $\mathcal{N}$. Tank $i$ has water inflow $Q_i^{\textrm{in}}$ and outflow $Q_i^{\textrm{out}}$. The concentrations of the substrate and biomass inflows are $S_i^{\textrm{in}}$ and $X_i^{\textrm{in}}$. The tank is assumed to be perfectly mixed and has substrate and biomass concentrations $S_i$ and $X_i$. We suppress the subscript to denote the vectors of these quantities in $\mathbb{R}^n$. $V\in\mathbb{R}^{n\times n}$ is a diagonal matrix in which $V_{ii}$ is the volume of tank $i$.

The substrate in tank $i$ is converted to biomass at the rate $r(S_i,X_i)/y=\mu(S_i,X_i)X_i/y$, where the growth rate, $\mu(s,x)$, is positive for $s>0$ and $x>0$. We refer to $\mu(s,x)x$ as the kinetics. The conversion occurs with yield $y>0$. We will focus on the following two growth rates.
\begin{itemize}
\item Monod~\cite{monod1949growth}:
\[
\mu_{\textrm{M}}(s)=\frac{\mu^{\max} s}{K+s}.
\]
\item Contois~\cite{contois1959kinetics}:
\[
\mu_{\textrm{C}}(s,x)=\frac{\mu^{\max}s}{Kx+s}.
\]
\end{itemize}
We let all tanks have the same $\mu^{\max}$, $K$, and $y$, and note that our results straightforwardly extend to the case where they are not identical due to, for instance, temperature or pH variation. Some of our results apply to other growth rates as well. We refer the reader to Appendix 1 of~\cite{bastin2013line} for a comprehensive list.

The Monod growth rate, $\mu_{\textrm{M}}(s)$, is concave. The corresponding kinetics, $\mu_{\textrm{M}}(s)x$, are quasiconcave, but not concave. The Contois growth rate, $\mu_{\textrm{C}}(s,x)$, is not concave. Observe that the corresponding kinetics can be written $\mu_{\textrm{C}}(s,x)x=\mu_{\textrm{M}}(s/x)x$. This is the perspective of the Monod growth rate. The kinetics corresponding to the Contois growth rate are therefore concave because the perspective of a concave function is always concave (see, e.g., Section 3.2.6 of~\cite{boyd2004convex}).

We let $Q_{ij}$ denote the flow from tank $i$ to tank $j$. Flow conservation implies that
\begin{equation}
Q_i^{\textrm{in}}+\sum_{j\in\mathcal{N}}Q_{ji}=Q_i^{\textrm{out}}+\sum_{j\in\mathcal{N}}Q_{ij}\label{MC}
\end{equation}
for $i\in\mathcal{N}$. Let $\tilde{d}_{ij}$ denote the diffusion between tanks $i$ and $j$, where $\tilde{d}_{ij}=\tilde{d}_{ji}$. Note that this could represent the sum of the diffusions in multiple pipes. For example, if $d_{ij}$ is the diffusion in a pipe with flow from $i$ to $j$, and $d_{ji}$ for another with flow from $j$ to $i$, then $\tilde{d}_{ij}=d_{ij}+d_{ji}$. We encounter this scenario in Section~\ref{disjunctive}, in which one has the decision to build a pipe in either direction.

Let $C=\textrm{diag}\left[Q_i^{\textrm{in}}\right]$, $G=\textrm{diag}\left[Q_i^{\textrm{out}}\right]$, and
\begin{align*}
M_{ij}&=\left\{
\begin{array}{ll}
Q_{ji}, & i\neq j\\
-Q_i^{\textrm{out}}-\sum_{k\in{\mathcal{N}}}Q_{ik},& i= j
\end{array}
\right.\\
&=\left\{
\begin{array}{ll}
Q_{ji},  & i\neq j\\
-Q_i^{\textrm{in}}-\sum_{k\in{\mathcal{N}}}Q_{ki},& i= j
\end{array}
\right.\\
L_{ij}&=\left\{
\begin{array}{ll}
\tilde{d}_{ij}, & i\neq j\\
-\sum_{k\in{\mathcal{N}}}\tilde{d}_{ik}, & i= j
\end{array}
\right..
\end{align*}
Let $\bm{1}\in\mathbb{R}^n$ be the vector of all ones. In matrix form (\ref{MC}) is given by $(M+C)\bm{1}=0$. Observe that because $L\bm{1}=0$, we also have $(M+L+C)\bm{1}=0$. Similarly, $\bm{1}^{\top}(M+G)=0$ and $\bm{1}^{\top}(M+L+G)=0$. 

The gradostat is a type of compartmental system, and $M$ is a compartmental matrix. A compartmental system is outflow connected if there is a directed path from every tank to a tank with outflow, i.e., a tank $i$ with $Q_i^{\textrm{out}}>0$. A key property of compartmental matrices is that the matrix $M$ is invertible if and only if the system is outflow connected. $M$ is irreducible if it cannot be made block lower triangular by reordering its indices. In the gradostat, this means that there is a directed path from each tank to every other tank. Note that if the gradostat is outflow connected and $M$ is irreducible, then $M+L$ is respectively invertible and irreducible. We refer the reader to~\cite{jacquez1993qualitative} for a thorough discussion of compartmental systems.

We henceforth assume that all gradostats are outflow connected. The matrices $M$ and $M+L$ are invertible, and therefore
\begin{subequations}
\begin{align}
-(M+L)^{-1}C\bm{1}&=\bm{1}\label{MLC}\\
-\left(M^{\top}+L\right)^{-1}G\bm{1}&=\bm{1}\label{MLG}.
\end{align}
\end{subequations}
We also know that $-M$ and $-M-L$ are $M$-matrices~\cite{jacquez1993qualitative}. This implies that $-M^{-1}$ and $-(M+L)^{-1}$ are nonnegative matrices. If $M$ is irreducible, then $-M^{-1}$ and $-(M+L)^{-1}$ are positive matrices~\cite{berman1994nonnegative}.

\subsection{Equilibria of the gradostat}\label{ss:gradostat}
The dynamics of the gradostat are
\begin{align*}
V\dot{S}&=-\frac{1}{y}VR(S,X) +(M+L)S + CS^{\textrm{in}}\\
V\dot{X}&= VR(S,X)+(M+L)X + CX^{\textrm{in}},
\end{align*}
where $R_i(S,X)=r(S_i,X_i)$, $i\in\mathcal{N}$. We obtain a steady state model by setting the derivatives to zero. We now briefly discuss when the solution to the steady state model corresponds to a unique, stable equilibrium of the gradostat.

Let $Z=X+yS$. The dynamics of $(Z,X)$ are
\begin{align*}
V\dot{Z}&=(M+L)Z + C\left(X^{\textrm{in}}+yS^{\textrm{in}}\right)\\
V\dot{X}&= VR((Z-X)/y,X)+(M+L)X + CX^{\textrm{in}}.
\end{align*}
Because the gradostat is outflow connected, $M+L$ is negative definite and $Z$ is globally asymptotically stable with equilibrium $\bar{Z}=-(M+L)^{-1}C\left(X^{\textrm{in}}+yS^{\textrm{in}}\right)$. Due to the cascade structure, the dynamics of $X$ are asymptotically equivalent to those obtained by replacing $Z$ with $\bar{Z}$. We assume that $\bar{Z}>0$; the following are two simple conditions that guarantee this.
\begin{itemize}
\item If $C\left(X^{\textrm{in}}+yS^{\textrm{in}}\right)$ is not the zero vector and $M$ is irreducible, then $(M+L)^{-1}$ is strictly negative~\cite{berman1994nonnegative} and $\bar{Z}>0$.
\item Because $M$ is a negative definite $M$-matrix, $(M+L)^{-1}$ is negative definite and nonpositive~\cite{berman1994nonnegative}. If $C\left(X^{\textrm{in}}+yS^{\textrm{in}}\right)>0$, i.e., there is inflow of substrate and/or biomass at every tank, then $\bar{Z}>0$.
\end{itemize}

We now briefly describe the equilibria of the gradostat, and refer the reader to \cite{rapaport2019notes} and Chapter 9 of \cite{smith1995theory} for more thorough discussions. If $X^{\textrm{in}}=0$, a `washout' equilibrium with $X=0$ always exists. A positive equilibrium with $S>0$ and $X>0$ exists if $\mu\left(\bar{Z}_i/y,0\right)>-M_{ii}-L_{ii}$ for all $i\in\mathcal{N}$. Intuitively, this means that the system can convert substrate to biomass faster than it ejects biomass. In this case, the washout equilibrium is repulsive. If $r\left(S_i,\bar{Z}_i-yS_i\right)$ is increasing and strictly concave on $\left[0,\bar{Z}_i/y\right]$, then there is at most one positive equilibrium. It is easy to verify that this last condition holds for Monod and Contois growth.

The steady state approximation by definition limits the range of scenarios we can consider. For example, a wastewater treatment system is clearly not in steady state during a storm surge. On the other hand, it may be a useful approximation when choosing where to install new pipes to improve efficiency under average operating conditions. Later in Section~\ref{EMP} we drop the steady state assumption and instead numerically approximate the derivatives. This enables us to optimize trajectories of the gradostat under time-varying conditions.

\subsection{Objectives}
The objective is to maximize the production of biogas. This amounts to maximizing the conversion of substrate to biomass at a subset of output tanks, $\mathcal{M}\subseteq\mathcal{N}$. The corresponding objective is
\begin{equation}
\max\;\sum_{i\in\mathcal{M}}V_{ii}r(S_i,X_i).\label{maxgrowth}
\end{equation}
All formulations in this paper are valid with the objective (\ref{maxgrowth}), but the theoretical results of Section~\ref{relax} might not hold if $\mathcal{M}\subset\mathcal{N}$. If $\mathcal{M}=\mathcal{N}$, mass conservation implies that
\[
\sum_{i\in\mathcal{N}}Q_i^{\textrm{in}}S_i^{\textrm{in}}=\sum_{i\in\mathcal{N}}Q_i^{\textrm{out}}S_i+\frac{1}{y}V_{ii}r(S_i,X_i).
\]
If the left hand side is fixed and $\mathcal{M}=\mathcal{N}$, (\ref{maxgrowth}) is equivalent to
\[
\min\;\sum_{i\in\mathcal{N}}Q_i^{\textrm{out}}S_i.
\]
This corresponds to minimizing the substrate leaving the network.

We hereon focus on maximizing the generic objective $\mathcal{F}(T)$, where $T=R(S,X)$. This corresponds to biogas production if $\mathcal{F}(T)=\sum_{i\in\mathcal{M}}V_{ii}T_i.$ It can also accommodate additional features, e.g., concavity could reflect diminishing returns due to limited ability to store biogas.

\subsection{Problem statement}
The full problem is given by
\begin{subequations}
\begin{align}
\mathcal{P}:\quad\max\quad&\mathcal{F}(T) \label{P0}\\
\textrm{such that}\quad&T_i=r(S_i,X_i),\quad i\in\mathcal{N}\label{P3}\\
&\frac{1}{y}VT= (M+L)S + CS^{\textrm{in}}\label{P1}\\
&-VT= (M+L)X + CX^{\textrm{in}}\label{P2}\\
&0=(M+C)\bm{1}\label{P4}\\
&\left(d,Q,Q^{\textrm{in}},S,S^{\textrm{in}},X,X^{\textrm{in}}\right)\in\Omega.\label{QUps}
\end{align}
\end{subequations}
The variables are $d$, $Q$, $Q^{\textrm{in}}$, $S$, $S^{\textrm{in}}$, $X$, $X^{\textrm{in}}$, and $T$. The growth constraint, (\ref{P3}), equates $T_i$ to the kinetics in tank $i$. (\ref{P1}), (\ref{P2}), and (\ref{P4}) balance the substrate, biomass, and water in each tank. The set $\Omega$ in (\ref{QUps}) represents generic constraints on design and/or operation. Several possibilities are listed in Use Case~\ref{ex:Omega} below. Note that the tank volumes are not variables, but could in principle be incorporated in a tractable manner using the techniques of Section~\ref{disjunctive}.

\begin{uc}\label{ex:Omega}
The following are operational constraints that could be represented by the set $\Omega$.
\begin{itemize}
\item $X^{\textrm{in}}=0$. All biomass in the system either is converted from substrate or was already present before the system came to steady state.
\item For each $i\in\mathcal{N}$, $Q^{\textrm{out}}_iS_i\leq\check{S}_i$. The mass of the substrate released from each tank per unit of time cannot exceed some limit.
\item $Q^{\textrm{in}\top}S^{\textrm{in}}=\bar{S}$. The total mass of substrate that enters the network per unit of time, $\bar{S}$, is allocated over the tanks. This could also be the case with biomass.
\item $\bm{1}^{\top}Q^{\textrm{in}}\leq\bar{Q}$. The total inflow (and, by conservation, ouflow) of water is limited. Observe that if the inflow is too small, the system will run inefficiently, but if it is too large, washout will be a stable equilibrium.
\item For each $i\in\mathcal{N}$, $\mu\left(\bar{Z}_i/y,0\right)>-M_{ii}-L_{ii}$. This ensures that the gradostat has a positive equilibrium and that the washout equilibrium, if it exists, is repulsive. Note that here $\bar{Z}$ is an optimization variable, subject to the constraints
\begin{align*}
0&=(M+L)\bar{Z}+C\left(X^{\textrm{in}}+yS^{\textrm{in}}\right)\\
\mu\left(\bar{Z}_i/y,0\right)&\geq\delta-M_{ii}-L_{ii},
\end{align*}
where $\delta$ is a small positive constant that makes the inequality strict. If $\mu(s,0)$ is concave in $s$, then this is a convex constraint. In the case of Contois growth, $\mu_{\textrm{C}}\left(\bar{Z}_i/y,0\right)$ is a constant. In the case of Monod growth, the latter constraint has an SOC representation, which has the same form as that given in~\ref{MonodX}.
\end{itemize}
\end{uc}

The growth constraint in $\mathcal{P}$, (\ref{P3}), is nonconvex. Sections~\ref{relax} and \ref{PRSOC} construct SOC relaxations of this constraint. We consider two cases for the rest of the problem.
\begin{itemize}
\item When the variables $d$, $Q$, and $Q^{\textrm{in}}$ are constant, (\ref{P1}) and (\ref{P2}) are linear. Then the resulting relaxations of $\mathcal{P}$ are SOCPs. In this case, the physical decisions are the concentrations of the substrate and biomass inflows, $S^{\textrm{in}}$ and $X^{\textrm{in}}$. Note that $S$, $X$, and $T$ are also optimization variables, but are not under direct control of a system operator.
\item When $d$, $Q$, and $Q^{\textrm{in}}$ are discrete variables, (\ref{P1}) and (\ref{P2}) are bilinear. We use disjunctive programming to linearize these constraints in Section~\ref{disjunctive}. In this case, the relaxations of $\mathcal{P}$ are MISOCPs. The physical decisions are $S^{\textrm{in}}$, $X^{\textrm{in}}$, the flows between tanks, $Q$, and the flows into the tanks, $Q^{\textrm{in}}$. $Q_{ij}$ could represent, for example, the decision to turn on a fixed speed pump in the pipe from tank $i$ to $j$, or the decision to build a pipe from $i$ to $j$. We give more detail on the forms of $d$, $Q$, and $Q^{\textrm{in}}$ in this case in Section~\ref{disjunctive}.
\end{itemize}

\begin{uc}\label{ex:soil}
Microbial activity in soil produces biogas, which, because it is not captured, contributes to the greenhouse effect. This is often modeled with Monod dynamics~\cite{schimel2003implications}. It has  also been shown to be dependent on spatial heterogeneity~\cite{parton1988dynamics,moorhead2006theoretical}, motivating the use of compartmental modeling. While there are no true design variables, it is of interest to understand which spatial structures lead to the greatest release of biogas. We can estimate this by maximizing a gradostat's biogas production over its water inputs, $Q^{\textrm{in}}$, and flows, $Q$ and $d$. The biomass in soil evolves slowly~\cite{moorhead2012theoretical} and is thus often treated as a constant parameter~\cite{parton1988dynamics,moorhead2006theoretical}. We model this by dropping constraint (\ref{P2}) from $\mathcal{P}$ and setting $X=X^{\textrm{c}}$. In Section~\ref{MonodX}, we show that in this case the kinetics corresponding to the Monod growth rate have an SOC representation.
\end{uc}

\section{Relaxation of the growth constraint}\label{relax}
The growth constraint, (\ref{P3}), is a nonlinear equality and hence nonconvex. We relax (\ref{P3}) by replacing it with the inequality
\begin{equation}
T_i\leq r(S_i,X_i),\quad i\in\mathcal{N}.\label{P3R}
\end{equation}
We refer to $\mathcal{P}$ with (\ref{P3R}) instead of (\ref{P3}) as $\mathcal{P}_{\textrm{R}}$. Unlike the equality (\ref{P3}), (\ref{P3R}) is convex for some growth rates, and can sometimes be represented as an SOC constraint; this is the focus of Section~\ref{PRSOC}. In this section, we analyze when $\mathcal{P}_{\textrm{R}}$ is exact, which is to say has the same optimal solution as $\mathcal{P}$.

First, observe that if the optimal solution to $\mathcal{P}_{\textrm{R}}$ satisfies (\ref{P3R}) with equality for all $i\in\mathcal{N}$, then $\mathcal{P}_{\textrm{R}}$ is exact. Given a solution to $\mathcal{P}_{\textrm{R}}$, we can therefore determine if it is feasible and optimal for $\mathcal{P}$ by simply checking if (\ref{P3R}) binds.

\begin{definition}
Let $\mathcal{E}\subseteq\mathcal{N}$ be such that for each $i\in\mathcal{E}$, $Q_{ji}+d_{ji}=0$ for all $j\in\mathcal{N}$, i.e., it receives no flow from other tanks. A gradostat is fully fed if for each $i\in\mathcal{E}$, $Q^{\textrm{in}}_i>0$, $S^{\textrm{in}}_i>0$, and $X^{\textrm{in}}_i>0$.
\end{definition}

\begin{lemma}\label{Spos}
Suppose the gradostat is fully fed and outflow connected. Then $S>0$ and $X>0$.
\end{lemma}
\begin{pf}
Consider a path, $\mathcal{L}$, from tank $s$ to tank $t$. Assume that $S^{\textrm{in}}_s>0$, $X^{\textrm{in}}_s>0$, and $Q^{\textrm{out}}_t>0$, and that for each edge $ij\in\mathcal{L}$, $Q_{ij}+d_{ij}>0$. We proceed by induction on the path.

Observe that because flow must enter every tank and $S^{\textrm{in}}_s>0$ and $X^{\textrm{in}}_s>0$, we must have $S_s>0$ and $X_s>0$ for (\ref{P1}) and (\ref{P2}) to be feasible at tank $s$. Now suppose that $jk\in\mathcal{L}$ and that $S_j>0$ and $X_j>0$. Then $S_k>0$ and $X_k>0$ for (\ref{P1}) and (\ref{P2}) to be feasible at tank $k$. Therefore, by induction, $S_i>0$ and $X_i>0$ for each $i\in\mathcal{L}$.

Because we have assumed that the gradostat is outflow connected and fully fed, all tanks must lie on a path like $\mathcal{L}$. Therefore, $S>0$ and $X>0$.\qed
\end{pf}

In an outflow connected gradostat, there can only be zero substrate  and biomass at tanks which have no inflow of substrate or biomass and receive no flow or diffusion from other tanks. Such tanks cannot exist in a fullly fed gradostat. Observe that another sufficient condition for $S>0$ is for the graph of diffusive couplings to be connected, which is equivalent to $\textrm{rank}(L)=n-1$.

\begin{assumption}\label{derS}
$r(s,x)$ is concave, differentiable, and for all $s\geq0$ and $x\geq0$, one has
\[
\frac{1}{y}\frac{\partial r(s,x)}{\partial s}-\frac{\partial r(s,x)}{\partial x}\geq0.
\]
\end{assumption}
We will discuss Assumption~\ref{derS} for specific growth rates in Section~\ref{PRSOC}.

In Theorem~\ref{PRexact} below, $S$, $X$, and $T$ are the only variables in $\mathcal{P}_{\textrm{R}}$, the remaining variables in $\mathcal{P}$ are regarded to be constant, and we drop the operational and design constraints in (\ref{QUps}). We will discuss how the theorem extends to the general case after the proof.

\begin{theorem}\label{PRexact}
$\mathcal{P}_{\textrm{R}}$ is exact if $\mathcal{F}(T)$ is concave and differentiable, the gradostat is outflow connected, Assumption~\ref{derS} holds, and either of the following is true:
\begin{itemize}
\item $M$ is irreducible and $\left(M^{\top}+L\right)V^{-1}\nabla\mathcal{F}(T)\leq 0$ and is not uniformly zero for all $T\geq0$; or
\item the gradostat is fully fed and $\left(M^{\top}+L\right)V^{-1}\nabla\mathcal{F}(T)< 0$ for all $T\geq0$.
\end{itemize}
\end{theorem}

\begin{pf}
$\mathcal{P}_{\textrm{R}}$ is convex due to Assumption~\ref{derS}. If $\mathcal{P}_{\textrm{R}}$ satisfies a constraint qualification such as Slater's condition, any optimal solution must satisfy the Karush-Kuhn-Tucker (KKT) conditions~\cite{boyd2004convex}.

$\mathcal{P}_{\textrm{R}}$ satisfies Slater's condition if there is a feasible solution for which (\ref{P3R}) is strict. Because the gradostat is outflow connected and either is fully fed or has an irreducible $M$, if there is any inflow of substrate and biomass, then we must have $S>0$ and $X>0$ for (\ref{P1}) and (\ref{P2}) to be feasible. Therefore $r(S_i,X_i)>0$ for all $i\in\mathcal{N}$. We obtain a feasible solution for which (\ref{P3R}) is strict by setting $T=0$. Therefore, a Slater point exists, and any optimal solution of $\mathcal{P}_{\textrm{R}}$ satisfies the KKT conditions.

Let $\rho\in\mathbb{R}^n_+$ be the vector of dual multipliers of constraint (\ref{P3R}), and let $\sigma\in\mathbb{R}^n$ and $\epsilon\in\mathbb{R}^n$ be the respective multipliers of (\ref{P1}) and (\ref{P2}). Let $U^S$ and $U^X$ be diagonal matrices with
\[
U^S_{ii}=\frac{\partial r(S_i,X_i)}{\partial S_i},\quad U^X_{ii}=\frac{\partial r(S_i,X_i)}{\partial X_i}
\]
for each $i\in\mathcal{N}$. The KKT conditions for $\mathcal{P}_{\textrm{R}}$ are given by
\begin{subequations}
\begin{align}
&(\ref{P1}),(\ref{P2}),(\ref{P3R})\nonumber\\
\rho&=\nabla\mathcal{F}(T)+V\left(\frac{1}{y}\sigma-\epsilon\right)\label{KKT1}\\
U^S\rho&=\left(M^{\top}+L\right)\sigma\label{KKT2}\\
U^X\rho&=\left(M^{\top}+L\right)\epsilon\label{KKT3}\\
0&=\rho_i\left(T_i-r(S_i,X_i)\right),\quad i\in\mathcal{N}.\label{KKTCS}
\end{align}
\end{subequations}
The complementary slackness condition, (\ref{KKTCS}), implies that if $\rho>0$, then constraint (\ref{P3R}) binds for all $i\in\mathcal{N}$ and $\mathcal{P}_{\textrm{R}}$ is exact.

Let $U=\frac{1}{y}U^S-U^X$ and $W=\left(M^{\top}+L\right)V^{-1}-U$. $U$ is positive semidefinite due to Assumption~\ref{derS}. Because the gradostat is outflow connected, $\left(M^{\top}+L\right)V^{-1}$ is negative definite, and therefore so is $W$. Arithmetic with (\ref{KKT1})-(\ref{KKT3}) yields
\begin{equation}
\rho=W^{-1}\left(M^{\top}+L\right)V^{-1}\nabla\mathcal{F}(T),\label{Wrho}
\end{equation}
Because $-W$ is a positive definite $M$-matrix, $W^{-1}$ is nonpositive. We have $\rho>0$ in both of the following two cases.
\begin{itemize}
\item If $M$ is irreducible, $W^{-1}$ is strictly negative~\cite{berman1994nonnegative}. Therefore, if $\left(M^{\top}+L\right)V^{-1}\nabla\mathcal{F}(T)\leq0$ and is not uniformly zero for all $T\geq0$, then $\rho>0$.
\item If $\left(M^{\top}+L\right)V^{-1}\nabla\mathcal{F}(T)<0$ for all $T\geq0$, then $\rho>0$ because $W^{-1}$ is negative definite and nonpositive.
\end{itemize}
These are the two conditions we assumed in the theorem. Therefore, $\rho>0$, constraint (\ref{P3R}) is met with equality, and $\mathcal{P}_{\textrm{R}}$ is exact.\qed
\end{pf}

We now discuss Theorem~\ref{PRexact} and its proof. We have regarded all variables except $S$, $X$, and $T$ to be fixed. The theorem therefore indicates when a solution to $\mathcal{P}_{\textrm{R}}$, $(S,X,T)$, solves equations (\ref{P3})-(\ref{P2}). Theorem~\ref{PRexact} also holds when $d$, $Q$, $Q^{\textrm{in}}$, $S^{\textrm{in}}$, and $X^{\textrm{in}}$ are variables, so long as its assumptions hold at their optimal values. If constraint (\ref{QUps}) only affects $d$, $Q$, $Q^{\textrm{in}}$, $S^{\textrm{in}}$, and $X^{\textrm{in}}$, it will not change the exactness. However, if (\ref{QUps}) constrains $S$ or $X$, then Theorem~\ref{PRexact} is not guaranteed to hold.



\begin{corollary}\label{coroF}
Suppose that 
\[
\mathcal{F}(T)=\sum_{i\in\mathcal{N}}V_{ii}T_i.
\]
$\mathcal{P}_{\textrm{R}}$ is exact if the gradostat is outflow connected, Assumption~\ref{derS} holds, and either of the following is true:
\begin{itemize}
\item $M$ is irreducible; or
\item the gradostat is fully fed and $Q^{\textrm{out}}> 0$.
\end{itemize}
\end{corollary}
\begin{pf}
The result is obtained by substituting $\nabla\mathcal{F}(T)=V\bm{1}$ in (\ref{Wrho}) and using (\ref{MLG}) to simplify the right hand side to $\rho=W^{-1}Q^{\textrm{out}}$.\qed
\end{pf}

There are likely further possible refinements and generalizations. For example, if $\mathcal{F}(T)$ is not differentiable, then the conditions of Theorem~\ref{PRexact} need to hold for all of its subgradients, or, alternatively, some subgradient at the optimal solution. On the other hand, there are surely systems of interest for which the relaxation is not exact. From this point of view, Theorem~\ref{PRexact} does not specify the entire set of gradostats for which $\mathcal{P}_{\textrm{R}}$ is an exact relaxation. Rather, it is theoretical evidence that $\mathcal{P}_{\textrm{R}}$ is exact for a meaningful class of gradostats and may be a high quality approximation for others.

\section{SOC representations of $\mathcal{P}_{\textrm{R}}$}\label{PRSOC}
In this section we construct original SOC representations and approximations of the growth constraint in $\mathcal{P}_{\textrm{R}}$, (\ref{P3R}). 

\subsection{Contois growth}
With the Contois growth rate, constraint (\ref{P3R}) takes the form
\begin{equation}
T_i\leq\frac{\mu^{\max}S_iX_i}{KX_i+S_i},\quad i\in\mathcal{N}.\label{ConsoisConcave}
\end{equation}
This is a convex constraint because, as discussed in Section~\ref{s:def}, the right hand side is concave.
\begin{theorem}\label{th:ContoisSOC}
(\ref{ConsoisConcave}) is equivalent to the hyperbolic constraint
\begin{subequations}\label{ContoisSOC}
\begin{align}
\tilde{S}_i^2+\tilde{T}_i^2&\leq \left(\tilde{S}_i-\tilde{T}_i\right)\left(2\mu^{\max}KX_i + \tilde{S}_i-\tilde{T}_i\right)\label{ContoisSOC1}\\
0&\leq \tilde{S}_i-\tilde{T}_i\label{ContoisSOC2},
\end{align}
\end{subequations}
where $\tilde{S}_i=\mu^{\max}S_i$ and $\tilde{T}_i=KT_i$, $i\in\mathcal{N}$.
\end{theorem}
This can be shown straightforwardly by simplifying (\ref{ContoisSOC1}). Note that (\ref{ContoisSOC2}) ensures the nonnegativity of both multiplicative terms on the right hand side of (\ref{ContoisSOC1}), and is implied by (\ref{ConsoisConcave}). We refer to $\mathcal{P}_{\textrm{R}}$ with (\ref{ContoisSOC}) in place of (\ref{P3R}) as $\mathcal{P}_{\textrm{RC}}$.

Hyperbolic constraints are a special case of SOC constraints. In standard SOC form, (\ref{ContoisSOC1}) is written
\[
\left\|\left[\begin{array}{c}
\tilde{S}_i\\
\tilde{T}_i\\
\mu_{\max}KX_i 
\end{array}\right]\right\|\leq \mu_{\max}KX_i +  \tilde{S}_i - \tilde{T}_i.
\]

We now test Assumption~\ref{derS} for $\mathcal{P}_{\textrm{RC}}$. We have
\[
\frac{1}{y}\frac{\partial r(s,x)}{\partial s}-\frac{\partial r(s,x)}{\partial x}=
\frac{\mu^{\max}\left(Kx^2-ys^2\right)}{y\left(Kx+s\right)^2}.
\]
Clearly this is negative for some $s\geq 0$ and $x\geq0$, and hence does not perfectly satisfy Assumption~\ref{derS}. However, given that there is usually more biomass than substrate, $y<1$, and typically $K\approx1$, we expect it to hold around the optimal solution for most realistic systems. Therefore, we expect Theorem~\ref{PRexact} to hold most of the time for $\mathcal{P}_{\textrm{RC}}$.

\subsection{Monod growth}
With the Monod growth rate, constraint (\ref{P3R}) takes the form
\begin{equation}
T_i\leq\frac{\mu^{\max}S_iX_i}{K+S_i},\quad i\in\mathcal{N}.\label{MonodQC}
\end{equation}
We refer to $\mathcal{P}_{\textrm{R}}$ with (\ref{MonodQC}) in place of (\ref{P3R}) as $\mathcal{P}_{\textrm{RM}}$. 

Because the right side of (\ref{MonodQC}) is not concave, Assumption~\ref{derS} does not fully hold, and Theorem~\ref{PRexact} does not directly apply. However, it is differentiable, which means that the KKT conditions hold at any local optimum. If the derivative condition in Assumption~\ref{derS} is satisfied along with the other conditions of Theorem~\ref{PRexact} or Corollary~\ref{coroF}, then $\mathcal{P}_{\textrm{RM}}$ is exact. Evaluating the derivative condition in Assumption~\ref{derS} yields
\[
\frac{1}{y}\frac{\partial r(s,x)}{\partial s}-\frac{\partial r(s,x)}{\partial x}=
\frac{\mu^{\max}\left(Kx^2-ys(K+s)\right)}{y\left(K+s\right)^2}.
\]
As with the Contois growth rate, this is negative for some $s\geq0$ and $x\geq0$, but is likely to be nonnegative around the solution for a realistic system. Therefore, we expect $\mathcal{P}_{\textrm{RM}}$ to be exact, albeit at a possibly local optimum.

Unfortunately, without a more tractable representation of (\ref{MonodQC}), $\mathcal{P}_{\textrm{RM}}$ will be difficult to solve at larger scales. In Sections~\ref{MonodX} and~\ref{MonodCE}, we construct two different SOC approximations of $\mathcal{P}_{\textrm{RM}}$.

\subsubsection{Constant biomass}\label{MonodX}
In some applications, the biomass, $X$, does not change significantly relative to $S$. We now assume that biomass is constant, i.e., $X=X^{\textrm{c}}$. Given this assumption, (\ref{P3R}) takes the form
\begin{equation}
T_i\leq\frac{\mu^{\max}S_iX^{\textrm{c}}_i}{K+S_i},\quad i\in\mathcal{N}.\label{MonodConcave}
\end{equation}
This is a convex constraint because the right hand side, a Monod function, is concave. 
\begin{theorem}\label{th:MonodSOC}
(\ref{MonodConcave}) is equivalent to the hyperbolic constraint
\begin{subequations}\label{MonodSOC}
\begin{align}
\hat{S}_i^2+\hat{T}_i^2&\leq \left(\hat{S}_i-\hat{T}_i\right)\left(2\mu^{\max}KX^{\textrm{c}}_i + \hat{S}_i-\hat{T}_i\right)\label{MonodSOC1}\\
0&\leq \hat{S}_i-\hat{T}_i\label{MonodSOC2},
\end{align}
\end{subequations}
where $\hat{S}_i=\mu^{\max}S_iX^{\textrm{c}}_i$ and $\hat{T}_i=KT_i$, $i\in\mathcal{N}$.
\end{theorem}
As with Theorem~\ref{th:ContoisSOC}, this can be shown by simplifying. We refer to $\mathcal{P}_{\textrm{R}}$ with (\ref{MonodSOC}) in place of (\ref{P3R}) and without the biomass balance, (\ref{P2}), as $\mathcal{P}_{\textrm{RM}X}$.
In standard SOC form, (\ref{MonodSOC1}) is written
\[
\left\|\left[\begin{array}{c}
\hat{S}_i\\
\hat{T}_i\\
\mu_{\max}KX^{\textrm{c}}_i 
\end{array}\right]\right\|\leq \mu_{\max}KX^{\textrm{c}}_i +  \hat{S}_i - \hat{T}_i.
\]

Theorem~\ref{PRexact} does not directly apply to $\mathcal{P}_{\textrm{RM}X}$ because $X=X^{\textrm{c}}$. The proof may be adapted by simply disregarding constraint (\ref{P2}) in $\mathcal{P}_{\textrm{R}}$; we omit the details because they are straightforward. Evaluating Assumption~\ref{derS} for (\ref{MonodConcave}) yields
\[
\frac{1}{y}\frac{\partial r(s,x)}{\partial s}-\frac{\partial r(s,x)}{\partial x}=
\frac{\mu^{\max}X^{\textrm{c}}K}{y\left(K+s\right)^2},
\]
which is always nonnegative. Therefore, (a slightly modified version of) Theorem~\ref{PRexact} always holds for $\mathcal{P}_{\textrm{RM}X}$. Given the approximation $X=X^{\textrm{c}}$, $\mathcal{P}_{\textrm{RM}X}$ is exact.

\subsubsection{Convex envelopes}\label{MonodCE}
As stated above, the relaxed Monod constraint (\ref{MonodQC}) does not lead to a tractable optimization problem. In this section, instead of making a physical approximation, we construct an SOC approximation of the non-relaxed Monod constraint,
\begin{equation}
T_i=\frac{\mu^{\max}S_iX_i}{K+S_i},\quad i\in\mathcal{N}.\label{MonodExact}
\end{equation}
Note that (\ref{MonodExact}) is an equality, whereas (\ref{MonodQC}) is an inequality. We relax this constraint using the concave and convex envelopes in Sections 3.1 and 3.2 of~\cite{tawarmalani2001semidefinite}; see also \cite{tawarmalani2002convex,tawarmalani2013convexification}.

We can rewrite (\ref{MonodExact}) as $\mu^{\max}X_i = T_i + KT_i/S_i$, $i\in\mathcal{N}$. The only nonconvexity is due to the term $T_i/S_i$. To apply the convex envelopes, we need upper and lower bounds of the form $\underline{S}_i \leq S_i \leq \overline{S}_i$ and $\underline{X}_i \leq X_i \leq \overline{X}_i$ for each $i\in\mathcal{N}$. Given such bounds, we use (\ref{MonodExact}), and the fact that the right hand side is increasing, to obtain the following upper and lower bounds on $T_i$:
\begin{equation}\label{Tlim}
\overline{T}_i=\frac{\mu^{\max}\overline{S}_i\overline{X}_i}{K+\overline{S}_i},\quad\underline{T}_i=\frac{\mu^{\max}\underline{S}_i\underline{X}_i}{K+\underline{S}_i}.
\end{equation}
We know that $S\geq0$, $X\geq0$, and therefore $T\geq0$. We now derive several other bounds.

\begin{lemma}\label{SXbounds}
For each $i\in\mathcal{N}$, in steady state,
\begin{enumerate}
\item $S_i\leq \max_{j\in\mathcal{N}}S_j^{\textrm{in}}$,
\item $X_i\geq \min_{j\in\mathcal{N}}X_j^{\textrm{in}}$, and
\item $X_i\leq \max_{j\in\mathcal{N}}X_j^{\textrm{in}}+yS_j^{\textrm{in}}$.
\end{enumerate}
\end{lemma}

\begin{pf}
We proceed casewise.
\begin{enumerate}
\item Starting from (\ref{P1}), because $T\geq0$ we have $(M+L)S+CS^{\textrm{in}}\geq0$, which implies that $(M+L)S+C\bm{1}\max_{j\in\mathcal{N}}S_j^{\textrm{in}}\geq0$. Applying (\ref{MLC}), we have $S\leq\bm{1}\max_{j\in\mathcal{N}}S_j^{\textrm{in}}$, i.e., $S_i\leq \max_{j\in\mathcal{N}}S_j^{\textrm{in}}$ for $i\in\mathcal{N}$.

\item Starting from (\ref{P2}), because $T\geq0$ we have $(M+L)X+CX^{\textrm{in}}\leq0$, which implies that $(M+L)X+C\bm{1}\min_{j\in\mathcal{N}}X_j^{\textrm{in}}\leq0$. As above, this leads to $X\geq\bm{1}\min_{j\in\mathcal{N}}X_j^{\textrm{in}}$, i.e., $X_i\geq \min_{j\in\mathcal{N}}X_j^{\textrm{in}}$ for $i\in\mathcal{N}$.

\item We combine (\ref{P1}) and (\ref{P2}) to obtain $(M+L)Z+CZ^{\textrm{in}}=0$, where $Z=X+yS$ and $Z^{\textrm{in}}=X^{\textrm{in}}+yS^{\textrm{in}}$. From here, similar to the previous cases, we can show that for each $i$, $\min_{j\in\mathcal{N}} Z^{\textrm{in}}_j\leq X_i+yS_i\leq \max_{j\in\mathcal{N}}Z^{\textrm{in}}_j$. Because $S\geq0$, this implies that for all $i\in\mathcal{N}$, $X_i\leq \max_{j\in\mathcal{N}} X^{\textrm{in}}_j+yS^{\textrm{in}}_j$.\qed
\end{enumerate}
\end{pf}

Note that if $S^{\textrm{in}}$ and $X^{\textrm{in}}$ are variables, we can simply replace them with their upper and lower bounds in Lemma~\ref{SXbounds}. We can now apply the convex envelopes of~\cite{tawarmalani2001semidefinite}. For each $i\in\mathcal{N}$, the Monod equation is represented by
\begin{subequations}\label{ConvEnv}
\begin{equation}
\mu^{\max}X_i= T_i + K\beta_i.
\end{equation}
$\beta_i$ represents the term $T_i/S_i$. It is constrained by the convex envelope, given below. The concave overestimator is the pair of linear inequalities
\begin{align}
\beta_i\underline{S}\overline{S}&\leq\overline{S}T_i-\underline{T}S_i+\underline{S}\underline{T}\\
\beta_i\underline{S}\overline{S}&\leq\underline{S}T_i-\overline{T}S_i+\overline{S}\overline{T}.
\end{align}
The convex underestimator is given by
\begin{align}
&\gamma_i\psi_i\geq \underline{T}\left( \frac{\overline{T}-T_i}{\overline{T}-\underline{T}} \right)^2\label{CEa}\\
&(\beta_i-\gamma_i)(S_i-\psi_i)\geq\overline{T}\left( \frac{T_i-\underline{T}}{\overline{T}-\underline{T}} \right)^2\\
&\psi_i\geq \max\left\{\underline{S}\frac{\overline{T}-T_i}{\overline{T}-\underline{T}}, S_i-\overline{S}\frac{T_i-\underline{T}}{\overline{T}-\underline{T}} \right\}\\
&\psi_i\leq \min\left\{\overline{S}\frac{\overline{T}-T_i}{\overline{T}-\underline{T}}, S_i-\underline{S}\frac{T_i-\underline{T}}{\overline{T}-\underline{T}} \right\}\\
&\beta_i-\gamma_i\geq0,\quad\gamma_i\geq0.
\end{align}
\end{subequations}
The first two constraints are hyperbolic SOC like (\ref{ContoisSOC}) and (\ref{MonodSOC}), and the rest are linear. $\psi_i$ and $\gamma_i$ are auxiliary variables.

We refer to $\mathcal{P}$ with (\ref{ConvEnv}) instead of (\ref{P3}) as $\mathcal{P}_{\textrm{RME}}$. $\mathcal{P}_{\textrm{RME}}$ may be a very good approximation to $\mathcal{P}$, but there are no theoretical results guaranteeing exactness.

\section{Linearization of bilinear terms}\label{disjunctive}
When $d$, $Q$, and $Q^{\textrm{in}}$ are variables, the mass flow terms in constraints (\ref{P1}) and (\ref{P2}) are bilinear. Bilinear constraints are nonconvex and in general difficult to optimize over. Two common ways to deal with bilinearities are convex relaxations, e.g., McCormick~\cite{mccormick1976computability} and lift-and-project~\cite{Sherali1992bilinear}, and disjunctive programming~\cite{grossmann2002disj}. Disjunctive programming techniques have been used in several mathematically similar problems, including chemical process optimization~\cite{lee2003global} and transmission network expansion in power systems~\cite{Pereira2001mixed}. Here we let the flows be discrete, in which case disjunctive programming leads to exact linearization. We note that relaxations are more appropriate when the flows are continuous instead of discrete.
%

We denote the set of pipes $\mathcal{J}$. For each pipe $ij\in\mathcal{J}$, define the binary variable
\begin{equation}
\lambda_{ij}\in\{0,1\}.\label{lambdabin}
\end{equation}
This could represent the decision to build a pipe, or the decision to turn on a fixed-speed pump. If both $ij\in\mathcal{J}$ and $ji\in\mathcal{J}$, we prohibit simultaneous flow in both directions with the constraint
\begin{equation}
\lambda_{ij}+\lambda_{ji}\leq 1.\label{flow1way}
\end{equation}

The flow rate through pipe $ij\in\mathcal{J}$ is given by $Q_{ij}=Q_{ij}^0 + \lambda_{ij}Q_{ij}^1$, where the constants $Q_{ij}^0$ and $Q_{ij}^1$ are respectively the base flow and the added flow if $\lambda_{ij}=1$. Introduce the variable $F_{ij}^S$ and the bilinear constraint $F_{ij}^S=\lambda_{ij}Q_{ij}^1S_i$. The flow of substrate from $i$ to $j$ is $Q_{ij}^0S_i+F_{ij}^S$. Because $\lambda_{ij}\in\{0,1\}$, we can rewrite this as the pair of linear disjunctive \cite{grossmann2002disj} constraints 
\begin{subequations}\label{Fdisj}
\begin{equation}
(1-\lambda_{ij})\Gamma\geq\left|Q_{ij}^1S_i - F_{ij}^S\right|\label{DJS1},\quad\lambda_{ij}\Gamma\geq\left|F_{ij}^S\right|,
\end{equation}
where $\Gamma$ is a large positive number. Hence $F_{ij}^S=0$ when $\lambda_{ij}=0$, and $F_{ij}^S=Q_{ij}^1S_i$ when $\lambda_{ij}=1$. We can similarly represent the constraint $F_{ij}^X=\lambda_{ij}Q_{ij}^1X_i$ as the pair of linear constraints
\begin{equation}
(1-\lambda_{ij})\Gamma\geq\left|Q_{ij}^1X_i - F_{ij}^X\right|,\quad \lambda_{ij}\Gamma\geq\left|F_{ij}^X\right|.\label{DJX2}
\end{equation}
\end{subequations}

The diffusion of substrate in pipe $ij\in\mathcal{J}$ is $\left(d_{ij}^0 + \lambda_{ij}d_{ij}^1\right)(S_i-S_j)$. Introduce the variables $G^S_{ij}$ and $G^X_{ij}$ and the bilinear constraints $G_{ij}^S=\lambda_{ij}d_{ij}^1(S_i-S_j)$ and $G_{ij}^X=\lambda_{ij}d_{ij}^1(X_i-X_j)$. We can similarly represent these as
\begin{subequations}\label{Gdisj}
\begin{align}
&(1-\lambda_{ij})\Gamma\geq\left|d_{ij}^1(S_i-S_j) - G_{ij}^S\right|,\;\lambda_{ij}\Gamma\geq\left|G_{ij}^S\right|\label{DdS2}\\
&(1-\lambda_{ij})\Gamma\geq\left|d_{ij}^1(X_i-X_j) - G_{ij}^X\right|,\;\lambda_{ij}\Gamma\geq\left|G_{ij}^X\right|.\label{DdX2}
\end{align}
\end{subequations}

If $S^{\textrm{in}}_i$ and $X^{\textrm{in}}_i$ are variables, then the products $Q^{\textrm{in}}_iS^{\textrm{in}}_i$ and $Q^{\textrm{in}}_iX^{\textrm{in}}_i$ are bilinear as well. We linearize them by noting that, due to flow conservation, (\ref{MC}), and the fact that $Q^{\textrm{out}}$ is fixed, $\lambda$ fully determines $Q^{\textrm{in}}$. Let $H_{ij}^{S,1}=\lambda_{ij}Q_{ij}^1S^{\textrm{in}}_i$, $H_{ij}^{S,2}=\lambda_{ji}Q_{ji}^1S^{\textrm{in}}_i$, $H_{ij}^{X,1}=\lambda_{ij}Q_{ij}^1X^{\textrm{in}}_i$, and $H_{ij}^{X,2}=\lambda_{ji}Q_{ji}^1X^{\textrm{in}}_i$. These are equivalent to
\begin{subequations}\label{Hdisj}
\begin{align}
(1-\lambda_{ij})\Gamma&\geq \left| Q_{ij}^1S^{\textrm{in}}_i-H_{ij}^{S,1}\right|,\;\lambda_{ij}\Gamma\geq \left|H_{ij}^{S,1}\right|\\
(1-\lambda_{ji})\Gamma&\geq \left| Q_{ji}^1S^{\textrm{in}}_i-H_{ij}^{S,2}\right|,\;\lambda_{ji}\Gamma\geq \left|H_{ij}^{S,2}\right|\\
(1-\lambda_{ij})\Gamma&\geq \left| Q_{ij}^1X^{\textrm{in}}_i-H_{ij}^{X,1}\right|,\;\lambda_{ij}\Gamma\geq \left|H_{ij}^{X,1}\right|\\
(1-\lambda_{ji})\Gamma&\geq \left| Q_{ji}^1X^{\textrm{in}}_i-H_{ij}^{X,2}\right|,\;\lambda_{ji}\Gamma\geq \left|H_{ij}^{X,2}\right|.
\end{align}
\end{subequations}

After making the appropriate substitutions, we obtain the following replacements for constraints (\ref{P1})-(\ref{P4}), which we write in scalar form for each $i\in\mathcal{N}$. For clarity, we indicate beneath each new term the corresponding entry in the original constraint. The flow balance, (\ref{P4}), becomes
\begin{subequations}\label{disjMC}
\begin{equation}
Q_i^{\textrm{in}}+\sum_{j\in\mathcal{N}}\underset{Q_{ji}}{\underbrace{Q_{ji}^0+\lambda_{ji}Q_{ji}^1}}= Q_i^{\textrm{out}}+\sum_{j\in\mathcal{N}}\underset{Q_{ij}}{\underbrace{Q_{ij}^0+\lambda_{ij}Q_{ij}^1}}.
\end{equation}
The substrate balance, (\ref{P1}), becomes
\begin{align}
&\frac{1}{y}V_{ii}T_i+Q_{i}^{\textrm{out}}S_i+ \sum_{j\in\mathcal{N}}\underset{Q_{ij}S_i}{\underbrace{Q_{ij}^0S_i+F_{ij}^S}}+\underset{d_{ij}(S_i-S_j)}{\underbrace{d_{ij}^0(S_i-S_j) +G_{ij}^S}}  \nonumber\\
&\quad = \underset{Q^{\textrm{in}}_iS^{\textrm{in}}_i}{\underbrace{Q^{\textrm{out}}_iS^{\textrm{in}}_i+\sum_{j\in\mathcal{N}}\left(Q_{ij}^0-Q_{ji}^0\right)S^{\textrm{in}}_i+H_{ij}^{S,1}-H_{ij}^{S,2}}}\nonumber\\
&\quad+ \sum_{j\in\mathcal{N}}\underset{Q_{ji}S_j}{\underbrace{Q_{ji}^0S_j+F_{ji}^S}}+\underset{d_{ji}(S_j-S_i)}{\underbrace{d_{ji}^0(S_j-S_i)+G_{ji}^S}}.\label{P1scalar}
\end{align}
The biomass balance, (\ref{P2}) takes the same form as (\ref{P1scalar}), but with $X$ in place of $S$.
\end{subequations}

We now make several comments. $Q^{\textrm{out}}_i$ and $V_{ii}$ are constant. If they were instead binary variables, we could linearize the resulting bilinearities using the same technique. The same binary variables, $\lambda$, appear throughout (\ref{lambdabin})-(\ref{disjMC}). We could also straightforwardly generalize this to integer capacities by associating multiple binary variables with each pipe.

\section{Extensions}\label{extensions}
We now consider two basic extensions to $\mathcal{P}_{\textrm{R}}$ and its SOC representations.

\subsection{Convex underestimators of the growth constraint}\label{ConvexUE}
When $\mathcal{P}_{\textrm{R}}$ is not an exact relaxation, it is useful to constrain $T$ from below using an underestimator of the kinetics. To retain tractability, the underestimator should be a convex function. There are multiple ways to do this, such as using lift-and-project relaxations~\cite{Sherali1992bilinear} or convex envelopes~\cite{tawarmalani2002convex}.

In this section, we design simple underestimators using bounds of the form $\underline{S}\leq S\leq\overline{S}$, $\underline{X}\leq X\leq\overline{X}$, and $\underline{T}\leq T\leq\overline{T}$. These bounds could be from Lemma~\ref{SXbounds}, which makes no assumption about the growth rate except that it is nonnegative, and therefore applies for both Monod and Contois. Alternatively, such bounds might be implied by operational constraints in (\ref{QUps}).

We first consider the Contois growth rate, the kinetics of which are increasing in both $S_i$ and $X_i$. Let
\[
\underline{T}_i=\frac{\mu^{\max} \underline{S}_i\underline{X}_i}{K\underline{X}_i+\underline{S}_i},\;\overline{T}^S_i=\frac{\mu^{\max} \overline{S}_i\underline{X}_i}{K\underline{X}_i+\overline{S}_i},\;\overline{T}^X_i=\frac{\mu^{\max} \underline{S}_i\overline{X}_i}{K\overline{X}_i+\underline{S}_i}.
\]
For each $i\in\mathcal{N}$, we have
\begin{align}
T_i-\underline{T}_i\geq\max&\left\{ \frac{\overline{T}_i^S-\underline{T}_i}{\overline{S}_i-\underline{S}_i}\left(S_i-\underline{S}_i\right) ,\right.\nonumber\\
&\;\;\;\left.\frac{\overline{T}_i^X-\underline{T}_i}{\overline{X}_i-\underline{X}_i}\left(X_i-\underline{X}_i\right)\right\}.\label{OEContois}
\end{align}
The first argument of the maximum is a linear interpolator with $X$ fixed at its lower bound, and the latter with $S$ fixed at its lower bound.

Now consider the Monod growth rate. In the exact case, an underestimator is given by the convex envelope, (\ref{ConvEnv}). Under the constant biomass approximation in Section~\ref{MonodX}, $X=X^{\textrm{c}}$, the kinetics only depend on $S$, and we can again use the linear interpolator. Let $\underline{T}$ and $\overline{T}$ be as in (\ref{Tlim}). For each $i\in\mathcal{N}$, we have
\begin{equation}
T_i-\underline{T}_i\geq  \frac{\overline{T}_i-\underline{T}_i}{\overline{S}_i-\underline{S}_i}\left(S_i-\underline{S}_i\right).\label{OEMonod}
\end{equation}

\subsection{Multiple time periods}\label{EMP}
The steady state approximation in Section~\ref{ss:gradostat} is not appropriate if the system is undergoing a transient. Instead of setting the gradostat's derivatives to zero, we now replace them with numerical approximations that are linear in the variables.

Suppose that there are multiple time periods, indexed by the set $\mathcal{T}$ and each of length $\Delta$. We index each time-varying-quantity by $(t)$, e.g., $S(t)$ is the substrate concentration at time $t$. The objective is to maximize biogas over all time periods:
\begin{equation}
\sum_{t\in\mathcal{T}}\alpha^t\mathcal{F}(T(t)),\label{P0T}
\end{equation}
where $\alpha\in(0,1)$ is a discount factor.

Let $\mathcal{D}_t$ be a numerical approximation of the derivative at time $t$. For example, in the case of Euler's explicit method with time step $\Delta$, $\mathcal{D}_t[S]=(S(t+1)-S(t))/\Delta$. The substrate and biomass balances, (\ref{P1}) and (\ref{P2}), become
\begin{subequations}\label{PT}
\begin{align}
\mathcal{D}_t[S] &=-\frac{1}{y}VT(t) + (M+L)S(t) + CS^{\textrm{in}}(t)\\
\mathcal{D}_t[X] &=VT(t)+ (M+L)X(t) + CX^{\textrm{in}}(t)
\end{align}
\end{subequations}
for $t\in\mathcal{T}$. Note that there are more accurate choices for $\mathcal{D}_t$, e.g., Runge-Kutta schemes~\cite{betts1998survey}.

The network, as parametrized by $M$, $L$, and $C$, is constant, but could easily be made time-varying as well. The remaining constraints in $\mathcal{P}_{\textrm{R}}$ are simply enforced for all $t\in\mathcal{T}$.

Below are three scenarios modeled by this setup.
\begin{itemize}
\item $M$, $L$, and $C$ depend on a single vector of binary variables, $\eta$. This corresponds to designing the system, e.g., adding new pipes, so that its performance is optimized for a trajectory, e.g., a sequence of operating points in a representative day.
\item $M(t)$, $L(t)$, and $C(t)$ depend on a sequence of vectors of binary variables, $\eta(t)$, $t\in\mathcal{T}$. This corresponds to dynamically reconfiguring the system through time, e.g., choosing which valves to open or which fixed speed pumps to run in each time period.
\item $M$, $L$, and $C$ are constant. This corresponds to optimizing the trajectory of the substrate and biomass over time, subject to other operational constraints. The physical decisions, $S^{\textrm{in}}(t)$ and $X^{\textrm{in}}(t)$, $t\in\mathcal{T}$, represent schedules of substrate and biomass inflow concentrations.
\end{itemize}
It may be most natural to implement the last two scenarios via receding horizon control~\cite{mattingley2011receding}. In this case, only the decisions corresponding to the first time period are implemented. The time horizon is then pushed back by one period, the optimization is resolved, the `new' first period's decisions are implemented, and so on. This accommodates uncertainty by allowing the user to update the parameters, e.g., inflows and constraints on $S^{\textrm{in}}(t)$ and $X^{\textrm{in}}(t)$, as new information becomes available.

\section{Examples}\label{app}
We implement $\mathcal{P}_{\textrm{RC}}$, $\mathcal{P}_{\textrm{RM}X}$, and $\mathcal{P}_{\textrm{RME}}$ on numerical examples. We solve each optimization using the parser CVX~\cite{cvx} and the solver Gurobi~\cite{gurobi}. Table~\ref{tab:Ps} summarizes the features of each optimization model. The last column refers to whether or not the corresponding kinetics satisfy the derivative condition in Assumption~\ref{derS}.
\begin{table}[h]
	\label{tab:Ps}
	\centering
	\begin{tabular}{|l|l|l|l|}
	\hline
	Model & Growth rate & Class & Assump.~\ref{derS}  \\
	\hline
	$\mathcal{P}$ & Any & NLP & \\
	$\mathcal{P}_{\textrm{R}}$ & Any & NLP & \\
	$\mathcal{P}_{\textrm{RC}}$ & Contois & (MI)SOCP & Usually\\
	$\mathcal{P}_{\textrm{RM}}$ & Monod & NLP & Usually\\
	$\mathcal{P}_{\textrm{RM}X}$ & Monod, constant $X$ & (MI)SOCP & Always\\
	$\mathcal{P}_{\textrm{RME}}$ & Monod & (MI)SOCP & Unlikely\\
	\hline
	\end{tabular}
 	\caption{Summary of optimization models}
\end{table}

We measure the quality of the approximations in terms of the relative difference between the kinetics and the variable $T$:
\[
\mathcal{E}=\max_i \frac{\left|r(S_i,X_i)-T_i\right|}{r(S_i,X_i)}.
\]
If Theorem~\ref{PRexact} holds, then $\mathcal{E}$ will be zero for $\mathcal{P}_{\textrm{RC}}$ and/or $\mathcal{P}_{\textrm{RM}X}$. We expect $\mathcal{E}$ to usually be positive for $\mathcal{P}_{\textrm{RME}}$ because exactness is never guaranteed.

\subsection{Variable flows}
In these examples, the flows and diffusions depend on binary variables, as described in Section~\ref{disjunctive}. In this case, $\mathcal{P}_{\textrm{RC}}$, $\mathcal{P}_{\textrm{RM}X}$, and $\mathcal{P}_{\textrm{RME}}$ are MISOCPs. MISOCPs are NP-hard, but can be solved at moderate scales.

We first state the models in full. The objective in each case is to maximize biogas production, (\ref{maxgrowth}), with $\mathcal{M}=\mathcal{N}$. The constraints are as follows.
\begin{itemize}
\item A maximum budget, $\sum_{ij\in\mathcal{J}}c_{ij}\lambda_{ij}\leq B$,
where $c_{ij}$ is the cost to install a pipe from tank $i$ and $j$, and $B$ is the budget. This corresponds $\Omega$ in (\ref{QUps}).
\item For each $i\in\mathcal{N}$, the SOC growth constraint. If $\mathcal{P}_{\textrm{RC}}$, this is (\ref{ContoisSOC}); if $\mathcal{P}_{\textrm{RM}X}$, (\ref{MonodSOC}); if $\mathcal{P}_{\textrm{RME}}$, (\ref{ConvEnv}). In the case of $\mathcal{P}_{\textrm{RM}X}$, $X=X^{\textrm{c}}$, and the other constraints on $X$ are dropped.
\item For each $i\in\mathcal{N}$, a linear underestimator of the growth constraint. If $\mathcal{P}_{\textrm{RC}}$, this is (\ref{OEContois}); if $\mathcal{P}_{\textrm{RM}X}$, (\ref{OEMonod}).
\item Binary and linear disjuctive constraints on the flows and diffusions. For each $ij\in\mathcal{J}$, (\ref{lambdabin}), (\ref{flow1way}), (\ref{Fdisj}), (\ref{Gdisj}). Note that we do not include (\ref{Hdisj}) because $S^{\textrm{in}}$ and $X^{\textrm{in}}$ are constant in these examples.
\item Flow, substrate, and biomass balances. For each $i\in\mathcal{N}$, (\ref{disjMC}).
\end{itemize}

\subsubsection{Four tanks}
We first consider a small, four-tank example. The growth rate parameters (Monod and Contois) are $\mu_{\max}=K=y=1$. There is no base network, i.e., $Q^0=d^0=0$. All pairs of tanks are candidates for new pipes in either direction, so that there are twelve binary variables. For all $ij\in\mathcal{J}$, $Q^1_{ij}=1$ and $d^1_{ij}=0.3$. The cost of each new pipe is $c_{ij}=1$, and the budget is $B=4$. The tank parameters are $V=\textrm{diag}[1\;2\;3\;4]^{\top}$, $Q^{\textrm{out}}=[2\;1\;3\;2]^{\top}$, $S^{\textrm{in}}=[1\;3\;1\;2]^{\top}$, and $X^{\textrm{in}}=[4\;3\;2\;1]^{\top}$. In $\mathcal{P}_{\textrm{RM}X}$, we set $X^{\textrm{c}}=X^{\textrm{in}}$. We set $\Gamma=50$ in all disjunctive constraints.

The results are summarized in Table~\ref{tab:4tank}. $\mathcal{P}_{\textrm{RC}}$ and $\mathcal{P}_{\textrm{RM}X}$ are both exact, as predicted by Theorem~\ref{PRexact}, and both result in the same pipe additions. $\mathcal{P}_{\textrm{RME}}$ is not exact and has a slightly different solution.
\begin{table}[h]
	\centering
	\label{tab:4tank}
	\begin{tabular}{|l|l|l|l|l|}
	\hline
	Model & Time (s) & $\mathcal{E}$ & Objective & New pipes  \\
	\hline
	$\mathcal{P}_{\textrm{RC}}$ & 6.3  & 0 & 8.81 & 21, 23, 24, 43\\
	$\mathcal{P}_{\textrm{RM}X}$ & 2.8  & 0 & 10.21 & 21, 23, 24, 43\\
	$\mathcal{P}_{\textrm{RME}}$ & 5.3  & 2.2 & 15.87 & 21, 23, 24, 41\\
	\hline
	\end{tabular}
		\caption{Results for the four-tank system}
\end{table}

We now change the objective from $\sum_{i=1}^4V_{ii}T_i$ to $\sum_{i=2}^4V_{ii}T_i$, and leave all other parameters the same. In this case, Theorem~\ref{PRexact} does not apply. The results are summarized in Table~\ref{tab:4tankmod}.
\begin{table}[h]
	\label{tab:4tankmod}
	\centering
	\begin{tabular}{|l|l|l|l|l|}
	\hline
	Model & Time (s) & $\mathcal{E}$ & Objective & New pipes  \\
	\hline
	$\mathcal{P}_{\textrm{RC}}$ & 6.1  & 0.66 & 7.89 & 21, 23, 24, 43\\
	$\mathcal{P}_{\textrm{RM}X}$ & 2.8  & 0.49 & 8.55 & 21, 23, 24, 43\\
	$\mathcal{P}_{\textrm{RME}}$ & 5.2  & 2.15 & 14.62 & 21, 23, 24, 41\\
	\hline
	\end{tabular}
		\caption{Results for the four-tank gradostat with modified objective}
\end{table}
The pipe additions are unchanged, but now neither $\mathcal{P}_{\textrm{RC}}$ nor $\mathcal{P}_{\textrm{RM}X}$ are exact. In both, $T_1$ binds with its underestimator, (\ref{OEContois}) or (\ref{OEMonod}), and $T_i=r(S_i,X_i)$ for $i=2,3,4$. This indicates that the assumptions of Theorem~\ref{PRexact} are somewhat rigid, and that when they are violated, inexactness tends to occur locally.

\subsubsection{Wheel with $n$ tanks}
We now look at a larger example to see how the MISOCPs scale, using $\mathcal{P}_{\textrm{RC}}$ as the representative model.

The gradostat has $n$ tanks arranged in a wheel. The first tank is a central hub, and the other $n-1$ tanks are around rim. There is no base network. A pipe can be installed in either direction from the hub tank to any rim tank, i.e., $1i\in\mathcal{J}$ and $i1\in\mathcal{J}$ for $i\in\mathcal{N}\setminus 1$. A pipe can also be installed in either direction between each rim tank and its neighbor, i.e., $2n\in\mathcal{J}$, $n2\in\mathcal{J}$, and $i,i+1\in\mathcal{J}$ and $i+1,i\in\mathcal{J}$ for $i\in\mathcal{N}\setminus\{1,n\}$. If there are $n$ tanks, then there are $4n-1$ binary variables.

We consider an easy case and a hard case. In the easy case, for $i\in\mathcal{N}$, the volumes are $V_{ii}=i$ and the outflows $Q^{\textrm{out}}_i=1$. In the hard case, they are $V_{ii}=1+(i\mod 6)$ and $Q^{\textrm{out}}_i=1+(i\mod 7)$. The other parameters, which are the same in both cases, are  $S^{\textrm{in}}_i=i$, $X^{\textrm{in}}_i=n-i+1$, and maximum budget $B=1.5n$, and the rest are the same as in the previous example.

\begin{figure}[h]
 \centering
\includegraphics[width=\columnwidth]{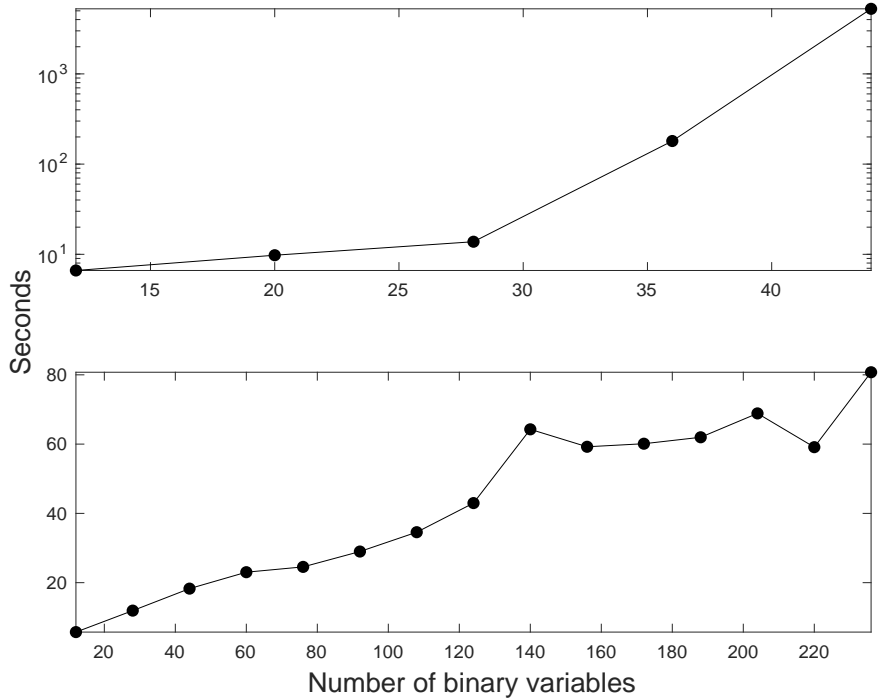} 
\caption{Computation time versus number of binary variables for $\mathcal{P}_{\textrm{RC}}$, an MISOCP, on a wheel gradostat. The hard case is shown on top and the easy case on bottom. Note that the top plot has a logarithmic $y$-axis.}
\label{fig:ExWheel}
\end{figure}

Figure~\ref{fig:ExWheel} show the results. In the easy case, the computation time increases roughly linearly, taking around a minute with with 60 tanks and 236 binary variables. In the hard case, the computation time increases exponentially, and more than an hour is needed with 11 tanks and 44 binary variables. 

The reason for the difference is that in the easy case, tanks with larger volumes are significantly better represented in the objective. This enables the solver to rapidly eliminate solutions with many pipes added to the smaller tanks, whereas in the harder case, the solver must search the feasible set more evenly. $\mathcal{P}_{\textrm{RC}}$ is exact in all cases.

\subsection{Constant flows and multiple time periods}
We now implement $\mathcal{P}_{\textrm{RC}}$ on an example with multiple time periods, as described in Section~\ref{EMP}. The water flows and diffusions are constant, which makes $\mathcal{P}_{\textrm{RC}}$ an SOCP. Our intention here is to demonstrate that this version of the problem can be solved at very large scales. This enables us to deal with larger systems and, as in this example, to choose the time step small enough that the continuous dynamics of the gradostat are well-represented.

There are four tanks, all with unit volume and growth rate parameters. The inflow vector is $Q^{\textrm{in}}=[2\;1\;1\;1]^{\top}$. The flows between tanks are: $Q_{12}=1$, $Q_{23}=2$, $Q_{34}=1$, $Q_{42}=1$, and the diffusion is $d=0.3Q$.

\begin{figure}[h]
 \centering
\includegraphics[width=\columnwidth]{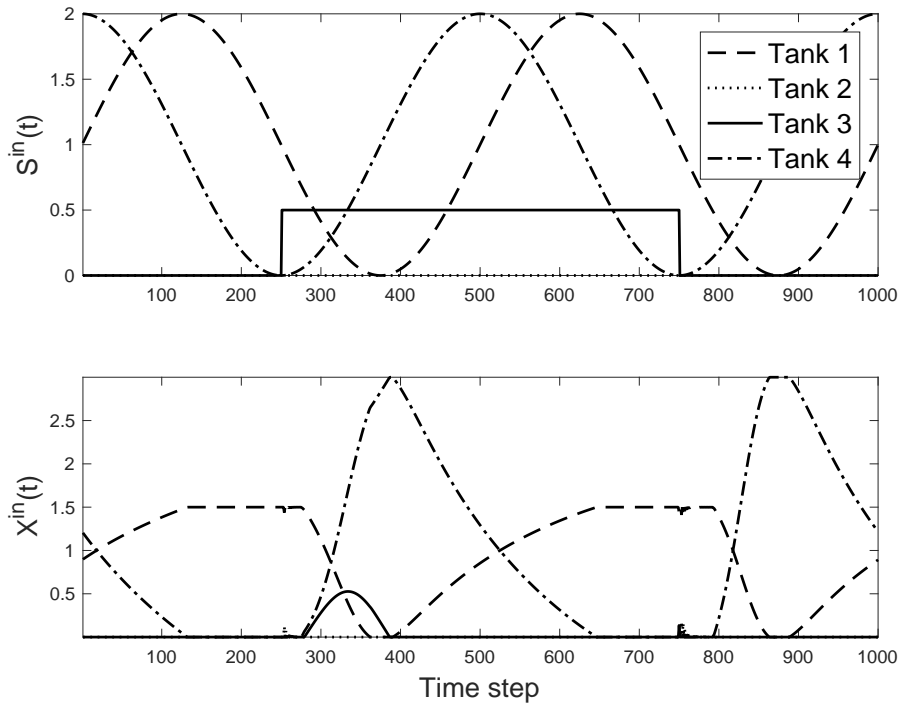} 
\caption{$S^{\textrm{in}}(t)$ and $X^{\textrm{in}}(t)$ as functions of time}
\label{fig:ExDynIn}
\end{figure}

There are $\tau=1000$ time periods of length $\Delta=1$. We approximate the derivative with Euler's method. The inflow substrate concentrations, shown in the top plot of Figure~\ref{fig:ExDynIn}, are $S^{\textrm{in}}_1(t)=1+\sin(4\pi t/\tau)$, $S^{\textrm{in}}_2(t)=0$,
\[
S^{\textrm{in}}_3(t)=\left\{
\begin{array}{cc}
1/2, & \tau/4< t \leq 3\tau/4\\
0, & \textrm{otherwise}
\end{array}
\right.,
\]
and $S^{\textrm{in}}_4(t)=1+\cos(4\pi t/\tau)$ for $t\in\mathcal{T}$.

The objective is to maximize the cumulative biogas production over all time periods, (\ref{P0T}). The constraints are listed below.
\begin{itemize}
\item At each time $t\in\mathcal{T}$, the total biomass added must satisfy $Q^{\textrm{in}\top}X(t)\leq 3$.
\item For each $i\in\mathcal{N}$ and $t\in\mathcal{T}$, the SOC Contois growth constraint, (\ref{ContoisSOC}).
\item At each time $t\in\mathcal{T}$, the dynamic substrate and biomass balances, (\ref{PT}). These constraints couple the variables in consecutive time periods.
\item As boundary conditions we require that $S(1)=S(\tau+1)$ and $X(1)=X(\tau+1)$.
\end{itemize}

There are roughly sixteen thousand variables in this problem. The computation time was 281 seconds, of which only 1.49 were taken by the solver and the rest by the parser. The optimal objective was 1140.18 units of biogas mass.

\begin{figure}[h]
 \centering
\includegraphics[width=\columnwidth]{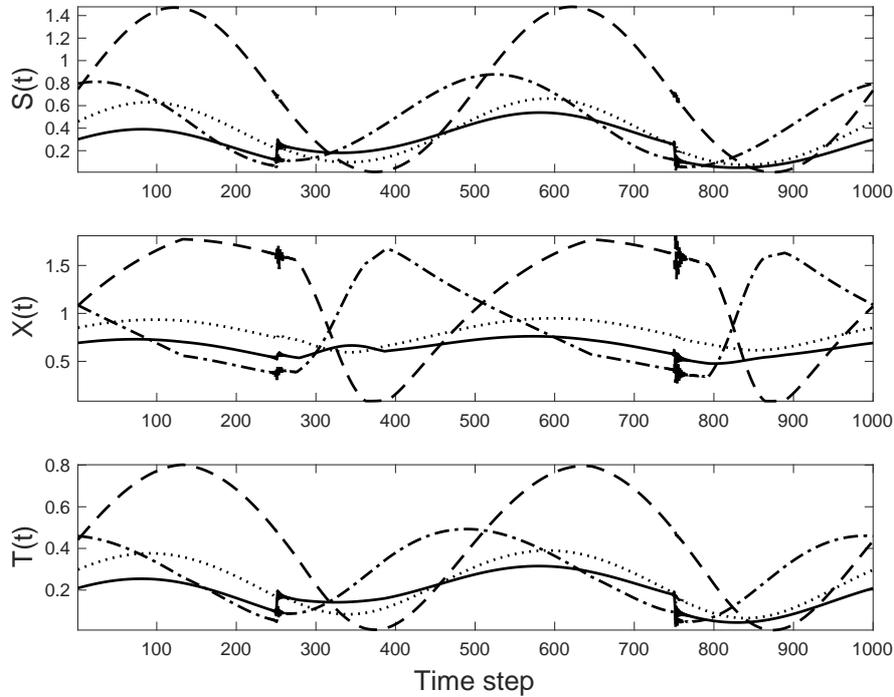} 
\caption{$S(t)$ and $X(t)$ as functions of time}
\label{fig:ExDynSXT}
\end{figure}

Figure~\ref{fig:ExDynIn} shows $S^{\textrm{in}}(t)$, which is as specified above, and the optimal $X^{\textrm{in}}(t)$. Figure~\ref{fig:ExDynSXT} shows the optimal $S(t)$, $X(t)$, and $T(t)$.

$X^{\textrm{in}}(t)$ is zero in tanks 2 and 3, except after $S^{\textrm{in}}_3(t)$ jumps up to $1/2$ and briefly when it drops back to zero. Notice that $\Delta$ is small enough for $X^{\textrm{in}}(t)$ and the curves in Figure~\ref{fig:ExDynSXT} to capture the oscillations caused by these transitions.

While $S^{\textrm{in}}_3(t)=1/2$, the substrate, biomass, and production in all tanks except the first jump, although slightly in tanks 2 and 4. This is because the additional substrate injected into tank 3 reaches tanks 2 and 4 through the flows. Tank 1 only receives substrate and biomass from tank 2 through a small diffusive coupling, and for this reason is not noticeably affected by the change in $S^{\textrm{in}}_3(t)$.

The solution was exact in all times periods. However, we observed that when there were fewer time periods and a larger time step, the solution was not always exact. We also remark that trajectories produced by $\mathcal{P}_{\textrm{RC}}$ do not necessarily lead to better stability or disturbance rejection. Such control objectives could be incorporated through a tracking objective and a receding horizon implementation. Characterizing the exactness of dynamic optimizations with different objectives is a topic of future work.

\section{Conclusions}
We have formulated SOCPs for optimizing the gradostat with Contois or Monod growth rates. The SOCPs are convex relaxations, which we proved are exact under simple conditions. We also gave linear underestimators, which are useful when the relaxations are not exact, and a dynamic extension in which the derivatives are replaced with numerical approximations instead of set to zero. 

There are many directions for future work. More physical features could be incorporated, including continuous variable flows via convex relaxations, recirculation of biomass, multi-reactions with several substrates, and other growth rates such as Teissier and Haldane. The constraints of the SOCPs could be used for a number of other purposes such as state estimation, setpoint tracking, and receding horizon control. As in this paper, the KKT conditions could be used to identify conditions under which each of these modifications is exact.

\section*{Acknowledgements}
We thank Professor Denis Dochain for helpful discussion.

\bibliography{MainBib,JATBib}

\begin{thebibliography}{10}
\expandafter\ifx\csname url\endcsname\relax
  \def\url#1{\texttt{#1}}\fi
\expandafter\ifx\csname urlprefix\endcsname\relax\def\urlprefix{URL }\fi
\expandafter\ifx\csname href\endcsname\relax
  \def\href#1#2{#2} \def\path#1{#1}\fi

\bibitem{shen2015overview}
Y.~Shen, J.~Linville, M.~Urgun-Demirtas, M.~Mintz, S.~Snyder, An overview of
  biogas production and utilization at full-scale wastewater treatment plants
  ({WWTP}s) in the {U}nited {S}tates: challenges and opportunities towards
  energy-neutral {WWTP}s, Renewable and Sustainable Energy Reviews 50 (2015)
  346--362.

\bibitem{monod1949growth}
J.~Monod, The growth of bacterial cultures, Annual review of microbiology 3~(1)
  (1949) 371--394.

\bibitem{contois1959kinetics}
D.~Contois, Kinetics of bacterial growth: relationship between population
  density and specific growth rate of continuous cultures, Microbiology 21~(1)
  (1959) 40--50.

\bibitem{Boyd1998SOCP}
M.~Lobo, L.~Vandenberghe, S.~Boyd, H.~Lebret, Applications of second-order cone
  programming, Linear Algebra and its Applications 284 (1998) 193--228.

\bibitem{tawarmalani2001semidefinite}
M.~Tawarmalani, N.~Sahinidis, Semidefinite relaxations of fractional programs
  via novel convexification techniques, Journal of Global Optimization 20~(2)
  (2001) 133--154.

\bibitem{drewesMISOCP2009}
S.~Drewes, Mixed integer second order cone programming, Ph.D. thesis,
  Technischen Universit\"{a}t Darmstadt, Department of Mathematics (2009).

\bibitem{Atamturk2010RC}
A.~Atamt\"{u}rk, V.~Narayanan, Conic mixed-integer rounding cuts, Mathematical
  Programming 122 (2010) 1--20.

\bibitem{belotti2017complete}
P.~Belotti, J.~G{\'o}ez, I.~P{\'o}lik, T.~Ralphs, T.~Terlaky, A complete
  characterization of disjunctive conic cuts for mixed integer second order
  cone optimization, Discrete optimization 24 (2017) 3--31.

\bibitem{gurobi}
{Gurobi Optimization, LLC}, \href{http://www.gurobi.com}{Gurobi optimizer
  reference manual} (2021).
\newline\urlprefix\url{http://www.gurobi.com}

\bibitem{smith1995theory}
H.~Smith, P.~Waltman, The theory of the chemostat: dynamics of microbial
  competition, Vol.~13, Cambridge University Press, 1995.

\bibitem{harmand2017chemostat}
J.~Harmand, C.~Lobry, A.~Rapaport, T.~Sari, The chemostat: Mathematical theory
  of microorganism cultures, John Wiley \& Sons, 2017.

\bibitem{lovitt1981gradostat}
R.~Lovitt, J.~Wimpenny, The gradostat: a bidirectional compound chemostat and
  its application in microbiological research, Microbiology 127~(2) (1981)
  261--268.

\bibitem{jager1987competition}
W.~J{\"a}ger, J.-H. So, B.~Tang, P.~Waltman, Competition in the gradostat,
  Journal of Mathematical Biology 25~(1) (1987) 23--42.

\bibitem{bayen2014minimal}
T.~Bayen, A.~Rapaport, M.~Sebbah, Minimal time control of the two tanks
  gradostat model under a cascade input constraint, {SIAM} Journal on Control
  and Optimization 52~(4) (2014) 2568--2594.

\bibitem{biegler1997systematic}
L.~Biegler, I.~Grossmann, A.~Westerberg, Systematic methods of chemical process
  design, Prentice Hall, 1997.

\bibitem{grossmann2002disj}
I.~Grossmann, Review of nonlinear mixed-integer and disjunctive programming
  techniques, Optimization and Engineering 3~(3) (2002) 227--252.

\bibitem{lee2003global}
S.~Lee, I.~Grossmann, Global optimization of nonlinear generalized disjunctive
  programming with bilinear equality constraints: applications to process
  networks, Computers \& Chemical Engineering 27~(11) (2003) 1557--1575.

\bibitem{tawarmalani2002convex}
M.~Tawarmalani, N.~Sahinidis, Convex extensions and envelopes of lower
  semi-continuous functions, Mathematical Programming 93~(2) (2002) 247--263.

\bibitem{tawarmalani2013convexification}
M.~Tawarmalani, N.~Sahinidis, Convexification and global optimization in
  continuous and mixed-integer nonlinear programming: theory, algorithms,
  software, and applications, Vol.~65, Springer Science \& Business Media,
  2013.

\bibitem{froment1990chemical}
G.~Froment, K.~Bischoff, J.~De~Wilde, Chemical reactor analysis and design, 3rd
  Edition, Wiley, 2010.

\bibitem{aris2000optimal}
R.~Aris, The Optimal Design of Chemical Reactors: A Study in Dynamic
  Programming, Academic Press, 1962.

\bibitem{bischoff1966optimal}
K.~Bischoff, Optimal continuous fermentation reactor design, The Canadian
  Journal of Chemical Engineering 44~(5) (1966) 281--284.

\bibitem{luyben1982optimal}
K.~C.~A. Luyben, J.~Tramper, Optimal design for continuous stirred tank
  reactors in series using {M}ichaelis--{M}enten kinetics, Biotechnology and
  Bioengineering 24~(5) (1982) 1217--1220.

\bibitem{hill1989minimum}
G.~Hill, C.~Robinson, Minimum tank volumes for {CFST} bioreactors in series,
  The Canadian Journal of Chemical Engineering 67~(5) (1989) 818--824.

\bibitem{de1996bioreactors}
C.~de~Gooijer, W.~Bakker, H.~Beeftink, J.~Tramper, Bioreactors in series: an
  overview of design procedures and practical applications, Enzyme and
  Microbial Technology 18~(3) (1996) 202--219.

\bibitem{harmand2003optimal}
J.~Harmand, A.~Rapaport, A.~Trofino, Optimal design of interconnected
  bioreactors: New results, AIChE Journal 49~(6) (2003) 1433--1450.

\bibitem{harmand2004optimal}
J.~Harmand, A.~Rapaport, A.~Dram{\'e}, Optimal design of two interconnected
  enzymatic bioreactors, Journal of Process Control 14~(7) (2004) 785--794.

\bibitem{harmand2005optimal}
J.~Harmand, D.~Dochain, The optimal design of two interconnected (bio) chemical
  reactors revisited, Computers \& chemical engineering 30~(1) (2005) 70--82.

\bibitem{zambrano2014optimizing}
J.~Zambrano, B.~Carlsson, Optimizing zone volumes in bioreactors described by
  {M}onod and {C}ontois growth kinetics, in: Proceeding of the IWA World Water
  Congress \& Exhibition, Lisbon, Portugal, 2014, p.~6.

\bibitem{zambrano2015optimal}
J.~Zambrano, B.~Carlsson, S.~Diehl, Optimal steady-state design of zone volumes
  of bioreactors with {M}onod growth kinetics, Biochemical engineering journal
  100 (2015) 59--66.

\bibitem{bayen2019steady}
T.~Bayen, P.~Gajardo, On the steady state optimization of the biogas production
  in a two-stage anaerobic digestion model, Journal of mathematical biology
  78~(4) (2019) 1067--1087.

\bibitem{crespo2020analysis}
M.~Crespo, A.~Rapaport, Analysis and optimization of the chemostat model with a
  lateral diffusive compartment, Journal of Optimization Theory and
  Applications 185 (2020) 597--621.

\bibitem{ROBLESRODRIGUEZ2018880}
C.~Robles-Rodriguez, J.~Bernier, V.~Rocher, D.~Dochain, A simple model of
  wastewater treatment plants for managing the quality of the {S}eine {R}iver,
  IFAC-PapersOnLine 51~(18) (2018) 880 -- 885, 10th IFAC Symposium on Advanced
  Control of Chemical Processes.
\newblock \href {https://doi.org/https://doi.org/10.1016/j.ifacol.2018.09.236}
  {\path{doi:https://doi.org/10.1016/j.ifacol.2018.09.236}}.

\bibitem{robles2019management}
C.~Robles-Rodriguez, A.~Ben-Ayed, J.~Bernier, V.~Rocher, D.~Dochain, Management
  of an integrated network of wastewater treatment plants for improving water
  quality in a river basin, in: IFAC Symposia Series, Vol.~52, 2019, p. 358.

\bibitem{bastin2013line}
G.~Bastin, D.~Dochain, On-line estimation and adaptive control of bioreactors,
  Vol.~1, Elsevier, 2013.

\bibitem{boyd2004convex}
S.~Boyd, L.~Vandenberghe, Convex Optimization, Cambridge University Press,
  2004.

\bibitem{jacquez1993qualitative}
J.~Jacquez, C.~Simon, Qualitative theory of compartmental systems, SIAM Review
  35~(1) (1993) 43--79.

\bibitem{berman1994nonnegative}
A.~Berman, R.~Plemmons, Nonnegative matrices in the mathematical sciences,
  SIAM, 1994.

\bibitem{rapaport2019notes}
A.~Rapaport, \href{https://hal.archives-ouvertes.fr/cel-02182956}{Lecture notes
  on mathematical models of interconnected chemostats} (Jul. 2019).
\newline\urlprefix\url{https://hal.archives-ouvertes.fr/cel-02182956}

\bibitem{schimel2003implications}
J.~Schimel, M.~Weintraub, The implications of exoenzyme activity on microbial
  carbon and nitrogen limitation in soil: a theoretical model, Soil Biology and
  Biochemistry 35~(4) (2003) 549--563.

\bibitem{parton1988dynamics}
W.~Parton, J.~Stewart, C.~Cole, Dynamics of {C}, {N}, {P} and {S} in grassland
  soils: a model, Biogeochemistry 5~(1) (1988) 109--131.

\bibitem{moorhead2006theoretical}
D.~Moorhead, R.~Sinsabaugh, A theoretical model of litter decay and microbial
  interaction, Ecological Monographs 76~(2) (2006) 151--174.

\bibitem{moorhead2012theoretical}
D.~Moorhead, G.~Lashermes, R.~Sinsabaugh, A theoretical model of {C}-and
  {N}-acquiring exoenzyme activities, which balances microbial demands during
  decomposition, Soil Biology and Biochemistry 53 (2012) 133--141.

\bibitem{mccormick1976computability}
G.~McCormick, Computability of global solutions to factorable nonconvex
  programs: {P}art {I}---{C}onvex underestimating problems, Mathematical
  programming 10~(1) (1976) 147--175.

\bibitem{Sherali1992bilinear}
H.~Sherali, A.~Alameddine, A new reformulation-linearization technique for
  bilinear programming problems, Journal of global optimization 2 (1992)
  379--410.

\bibitem{Pereira2001mixed}
L.~Bahiense, G.~Oliveira, M.~Pereira, A mixed integer disjunctive model for
  transmission network expansion, IEEE Transactions on Power Systems 16 (2001)
  560--565.

\bibitem{betts1998survey}
J.~Betts, Survey of numerical methods for trajectory optimization, Journal of
  Guidance, Control, and Dynamics 21~(2) (1998) 193--207.

\bibitem{mattingley2011receding}
J.~Mattingley, Y.~Wang, S.~Boyd, Receding horizon control, IEEE Control Systems
  Magazine 31~(3) (2011) 52--65.

\bibitem{cvx}
M.~Grant, S.~Boyd, {CVX}: Matlab software for disciplined convex programming,
  version 2.1, \url{http://cvxr.com/cvx} (Mar. 2014).

\end{thebibliography}

\end{document}